\documentclass[12pt]{amsart}
\makeatletter \@addtoreset{equation}{subsection} \makeatother

\usepackage{amssymb, euscript}
\usepackage{amsfonts, mathrsfs,}
\usepackage[matrix,arrow]{xy}

\newcommand{\lcm}{\operatorname{lcm}}
\newcommand{\discrep}{\operatorname{discrep}}

\newcommand{\rk}{\operatorname{rk}}
\newcommand{\Cl}{\operatorname{Cl}}

\newcommand{\Pic}{\operatorname{Pic}}
\newcommand{\Bs}{\operatorname{Bs}}

\newcommand{\comp}{\mathrel{\,\scriptstyle{\circ}\,}}

\newcommand{\qq}{\mathbin{\sim_{\scriptscriptstyle{\QQ}}}}

\newcommand{\down}[1]{\left\lfloor #1\right\rfloor}

\newcommand{\ov}[1]{\overline{#1}}

\newcommand{\OOO}{{\mathscr{O}}}
\newcommand{\MMM}{{\EuScript M}}
\newcommand{\EEE}{{\EuScript E}}
\newcommand{\FFF}{{\EuScript F}}

\newcommand{\GGG}{{\EuScript G}}
\newcommand{\HHH}{{\EuScript H}}
\newcommand{\B}{{\mathbf B}}
\newcommand{\PP}{\mathbb{P}}

\newcommand{\ZZ}{\mathbb{Z}}
\newcommand{\FF}{\mathbb{F}}

\newcommand{\CC}{\mathbb{C}}
\newcommand{\QQ}{\mathbb{Q}}

\newcommand{\xref}[1]{{\rm \ref{#1}}}

\newcounter{THN}[section]
\renewcommand{\theTHN}
{(\arabic{section}.\arabic{subsection})}
\newcounter{THNO}[section]
\renewcommand{\theTHNO}
{(\arabic{section}.\arabic{subsection}.\arabic{equation})}

\newenvironment{mparag}[1]{
\setcounter{THN}{\value{subsection}}
\refstepcounter{subsection}\refstepcounter{THN}
\par\medskip\noindent\begingroup \rm
{\bf\theTHN\ #1\ }}{\par\smallskip\endgroup}

\newenvironment{mtparag}[1]{
\setcounter{THN}{\value{subsection}}
\refstepcounter{subsection}\refstepcounter{THN}
\par\medskip\noindent\begingroup \it
{\bf\theTHN\ #1\ }}{\par\smallskip\endgroup}

\newenvironment{parag}[1]{
\setcounter{THNO}{\value{equation}}
\refstepcounter{equation}\refstepcounter{THNO}
\par\medskip\noindent\begingroup \rm
{\bf\theTHNO\ #1\ }}{\par\smallskip\endgroup}

\newenvironment{tparag}[1]{
\setcounter{THNO}{\value{equation}}
\refstepcounter{equation}\refstepcounter{THNO}
\par\medskip\noindent\begingroup \it
{\bf\theTHNO\ #1\ }}{\par\smallskip\endgroup}

\theoremstyle{definition}

\date{}

\author{Yu.~G.~Prokhorov}
\thanks{The author was 
partially supported by grants CRDF-RUM, No. 1-2692-MO-05 and
RFBR, No. 05-01-00353-a.}
\address{Department of Algebra, Faculty of Mathematics, Moscow State
University, Moscow 117234, Russia}
\email{prokhoro@mech.math.msu.su}
\title{The degree of $\QQ$-Fano threefolds}
\begin{document}

\maketitle

\section{Introduction}
\label{sect-1}
In this paper a \textit{$\QQ$-Fano variety} is a normal projective 
variety $X$ with at worst $\QQ$-factorial 
terminal singularities such that $-K_X$ is ample 
and $\Pic X$ is of rank one. Fano varieties with
terminal singularities form in 
important class because, according to the minimal 
model program, every variety of negative Kodaira dimension
should be birationally equivalent to a fibration $Y\to Z$
whose general fibre $Y_\eta$ belong to this class. Moreover,
in the case $\dim Z=0$, $Y_\eta=Y$ is of Picard number one,
i.e., $Y$ is a $\QQ$-Fano.

In dimension $2$ the only $\QQ$-Fano 
variety is the projective plane $\PP^2$.
In dimension $3$ $\QQ$-Fanos are bounded in the moduli sense 
by the following 
result of Kawamata: 
\begin{mtparag}{Theorem (\cite{Kawamata-1992bF}).} 
There exist positive integers $r$ and $d$ such that
for an arbitrary $\QQ$-Fano threefold $X$ 
we have $-K_X^3\le d$ and $rK_X$ is Cartier.
\end{mtparag}
Since the Weil divisor $-K_X$ gives a natural polarization of
a $\QQ$-Fano variety $X$, the rational number $-K_X^3$ is a very important invariant.
It is called \textit{the degree} of $X$.
In this paper we find a sharp bound for $-K_X^3$:
\begin{mtparag}{Theorem.}
\label{main-theorem}
Let $X$ be a $\QQ$-Fano threefold. 
Assume that $X$ is not Gorenstein.
Then $-K_X^3\le 125/2$ and the equality holds 
if and only if $X$ is isomorphic to the weighted projective space 
$\PP(1,1,1,2)$.
\end{mtparag}

Note that in the Gorenstein case we have the estimate
$-K_X^3\le 64$ by the classification of Iskovskikh and Mori-Mukai
and by Namikawa's result \cite{Namikawa-1997}.

The idea of the proof is as follows.
In Sections \ref{sect-4} and \ref{sect-5}
using Riemann-Roch formula for Weil divisors
\cite{Reid-YPG1987} and
Kawamata's estimates \cite{Kawamata-1992bF}
we produce a short list of possibilities
for singularities of $\QQ$-Fanos of degree
$\ge 125/2$. Here, to check a finite
(but very huge) number of Diophantine conditions,
we use a computer program
(cf. \cite{Suzuki-2004}). 
In Section \ref{sect-6} we exclude all these possibilities 
except for $\PP(1,1,1,2)$
by applying some birational transformations described in 
Section \ref{sect-3}. The techniques used 
on this step is a very common in birational geometry
(see \cite{Alexeev-1994ge}, \cite{Takagi-2002-I-II}, \cite{Takagi-2006}).
It goes back to Fano-Iskovskikh ``double projection method''.
The present paper is a logical continuation of 
our previous papers \cite{Prokhorov-2005a},
\cite{Prokhorov-2006-Enr} where we studied 
effective bounds of degree for sertain singular Fano 
threefolds.

\subsection*{Acknowledgements}
The work was carried out at Max-Planck-Institut f\"ur Mathematik,
Bonn in 2006. The author would like to thank the institute for the support
and hospitality. 


\section{Preliminaries}
\label{sect-2}
Throughout this paper, we work over the complex number field $\CC$.
\begin{mparag}{}
By $\Cl X$ we denote the Weil divisor class group 
of a normal variety $X$ (modulo linear equivalence).
There is a natural embedding $\Pic X \hookrightarrow \Cl X$.
Let $X$ be a Fano variety with at worst log terminal singularities.
It is well-known that both $\Pic X$ and $\Cl X$ 
are finitely generated and 
$\Pic X$ is torsion free (see e.g. \cite[\S 2.1]{Iskovskikh-Prokhorov-1999}). 
Moreover, numerical equivalence of
$\QQ$-Cartier divisors coincides with $\QQ$-linear one.
Therefore one can define the following numbers:
\[
\begin{array}{lll}
qF(X)&:=& \max \{ q \mid -K_X\qq q H,\quad H\in \Pic X\},
\\[8pt]
q\QQ(X)&:=& \max \{ q \mid -K_X\qq q L,\quad L\in \Cl X\},
\\[8pt]
qW(X)&:=& \max \{ q \mid -K_X\sim q L,\quad L\in \Cl X\}.
\end{array}
\]
By the above, all of them are positive, $q\QQ(X), qW(X)\in \ZZ$,
and $qF(X)\in \QQ$. If $X$ is smooth all these numbers 
coincide with the \textit{Fano index} of $X$. In general, we obviously have 
$q\QQ(X)\ge qF(X)$ and $q\QQ(X)\ge qW(X)$.

\begin{tparag}{Proposition (see e.g. {\cite[\S 2.1]{Iskovskikh-Prokhorov-1999}}).}
$qF(X)\le \dim X+1$.
\end{tparag}

The index $qW(X)$ was considered in \cite{Suzuki-2004}.
In particular, it was proved that $qW(X)\le 19$
for any $\QQ$-Fano threefold.
\end{mparag}


\begin{mparag}{Terminal singularities}
Let $(X,P)$ be a three-dimensional terminal singularity.
It follows from the classification that 
there is a one-parameter deformation 
$\mathfrak X\to \varDelta\ni 0$ over a small disk
$\varDelta\subset \CC$ such that the central fibre 
$\mathfrak X_0$ is isomorphic to $X$ and the generic fibre
$\mathfrak X_\lambda$ has only cyclic quotient singularities
$P_{\lambda,k}$ (see, e.g., \cite{Reid-YPG1987}). 
Thus, to every theefold $X$ with terminal singularities,
one can associate 
a collection $\B=\{(r_{P,k},b_{P,k})\}$, where
$P_{\lambda,k}\in \mathfrak X_\lambda$ is a singularity
of type $\frac 1{r_{P,k}}(b_{P,k},1,-1)$, $1\le b_{P,k}\le r_{P,k}/2$, 
$\gcd(r_{P,k}, b_{P,k})=1$.
This collection is uniquely determined by $X$  and called 
the \textit{basket} of singularities of $X$. By abuse of 
notation, we also will write $\B=(r_{P,k})$ instead of 
$\B=\{(r_{P,k},b_{P,k})\}$.
The index of $P$ is the least common multiple
of indices of points $P_{\lambda,k}$.
\end{mparag}

\begin{tparag}{Lemma (\cite[Corollary 5.2]{Kawamata-1988-crep}).}
\label{lemma-Kawamata-index-divisors}
Let $(X,P)$ be a three-dimensional terminal singularity
of index $r$ and let $D$ be a Weil $\QQ$-Cartier divisor on $X$.
There is an integer, $i$ such that $D\sim iK_X$ near $P$.
In particular, $rD$ is Cartier.
\end{tparag}
\begin{tparag}{Corollary.}
\label{cor-first-prop-index}
Let $X$ be a Fano threefold with terminal singularities
and let $r$ be the Gorenstein index of $X$. Then
\begin{enumerate}
\item 
$\gcd (r, qW(X))=1$,
\item 
$qF(X)r=q\QQ(X)$,
\item 
$qW(X)\le q\QQ(X)\le 4r$.
\end{enumerate}
\end{tparag}

\begin{parag}{}
Let $(X,P)$ be a three-dimensional terminal singularity
of index $r$ and let $D$ be a Weil $\QQ$-Cartier divisor on $X$.
By Lemma \ref{lemma-Kawamata-index-divisors}
there is an integer $i$ such that $0\le i<r$ and 
$D\sim iK_X$ near $P$. Deforming $D$ with $(X,P)$ 
we obtain Weil divisors $D_{\lambda}$ on $X_\lambda$.
Thus we have a collection of numbers $i_k$ such that
$0\le i_k<r_k$ and 
$D_{\lambda}\sim i_kK_{X_\lambda}$ near $P_{\lambda,k}$.
\end{parag}

\begin{mparag}{Riemann-Roch formula \cite{Reid-YPG1987}.}
Let $X$ be a threefold with terminal singularities 
and let $D$ be a Weil $\QQ$-Cartier divisor on $X$.
Then
\begin{multline}
\label{eq-RR}
\chi(D)=
\frac1{12}D\cdot (D-K_X)\cdot (2D-K_X)+
\\
+\frac1{12}D\cdot c_2+\sum_{P\in \B} c_P(D)+ 
\chi(\OOO_X),
\end{multline}
where 
\[
c_P(D)=-i_P\frac{r_P^2-1}{12r_P}+\sum_{j=1}^{i_P-1}\frac{\ov{b_Pj}(r_P-\ov{b_Pj})}{2r_P}.
\]
\end{mparag}

\begin{mparag}{}
\label{not-chi}
Now let $X$ be a Fano threefold with terminal singularities, 
let $q:=q\QQ(X)$, and
let $L$ be an ample Weil $\QQ$-Cartier divisor on $X$ such that $-K_X\qq qL$. 
By \eqref{eq-RR} we have
\begin{equation}
\label{eq-comput-chi}
\chi(tL)=1+\frac{t(q+t)(q+2t)}{12}L^3+
\frac{tL\cdot c_2}{12}
+ \sum_{P\in \B} c_P(tL),
\end{equation}
\[
c_P(tL)=
-i_{P,t} \frac{r_P^2-1}{12r_P}+\sum_{j=1}^{i_{P,t}-1}
\frac{\ov{b_Pj}(r_P-\ov{b_Pj})}{2r_P}.
\]
If $q>2$, then $\chi(-L)=0$. Using this equality
we obtain (see \cite{Suzuki-2004})
\begin{equation}
\label{eq-comput-L3}
L^3=\frac{12}{(q-1)(q-2)}\left( 
1-\frac{L\cdot c_2}{12}+\sum_{P\in B} c_P(-L)
\right).
\end{equation}
\end{mparag}

\begin{mparag}{}
In the above notation,
applying \eqref{eq-RR}, Serre duality and Kawamata-Viehweg vanishing 
to $D=K_X$ we get the following important equality (see, e.g., \cite{Reid-YPG1987}):
\begin{equation}
\label{eq-RR-O}
24=-K_X\cdot c_2+\sum_{P\in \B} \left( r_P-\frac 1{r_P}\right).
\end{equation}
Similarly, for $D=-K_X$ we have $H^i(X,-K_X)=0$ for $i>0$ and 
\[
c_P(-K_X)=\frac{r_P^2-1}{12r_P}-\frac{b_P(r-b_P)}{2r_P}.
\]
(see {\cite[\S 2]{Alexeev-1994ge}}). Combining this with \eqref{eq-RR-O} we obtain
\begin{equation}
\label{eq--K-exact}
\dim |-K_X|=
-\frac 12 K_X^3 +
2-\sum_{P\in \B} \frac{b_P(r_P-b_P)}{2r_P}. 
\end{equation}
In particular,
\begin{equation}
\label{eq--K-exact-1}
\dim |-K_X|\le 
-\frac 12 K_X^3 +
2-\frac 12 \sum_{P\in \B} \left(1- \frac1 {r_P}\right)\le 
-\frac 12 K_X^3 +
2.
\end{equation}
\begin{tparag}{Theorem ({\cite{Kawamata-1992bF}}, {\cite{KMMT-2000}}).}
\label{th-Kawamata-Kc}
In the above notation, $-K_X\cdot c_2(X)\ge 0$.
\end{tparag}
As a corollary we have ({\cite[\S 2]{Alexeev-1994ge}}):
\begin{equation}
\label{eq--K-exact-2}
\dim |-K_X|\ge 
-\frac 12 K_X^3 -2. 
\end{equation}

\begin{tparag}{Proposition ({\cite[\S 2]{Alexeev-1994ge}}]).}
\label{prop-A-lin-syst}
Let $X$ be a $\QQ$-Fano threefold. 
If $\dim |-K_X|\ge 2$, then the linear system 
$|-K_X|$ has no base components 
and is not composed of a pencil.
\textup(In particular, a general element of $|-K_X|$ 
is reduced and irreducible.\textup)
\end{tparag}
\end{mparag}

\begin{mparag}{}
\label{Kawamatas-ineq}
Now let $X$ be a $\QQ$-Fano threefold, 
let $q:=q\QQ(X)$, and
let $L$ be an ample Weil divisor on $X$ that generates 
the group $\Cl X/\mathrm{Tors}$. 
Let $\EEE$ be the double dual to $\Omega_X^1$. 
If $\EEE$ is not semistable, there is a maximal
destabilizing subsheaf $\FFF\subset \EEE$. 
Clearly, $c_1(\FFF)\equiv -pL$ for some $p\in \ZZ$.
Put $t:=p/q$, so that $c_1(\FFF)\equiv tK_X$.
According to 
\cite{Kawamata-1992bF} there are the following possibilities:
\begin{parag}{$\EEE$ is semistable.}
\label{mcase-semi}
Then $-K_X^3\le -3 K_X\cdot c_2(X)$.

\end{parag}

\begin{parag}{$\EEE$ is not semistable and $\rk \FFF=2$.}
\label{mcase-2}
Then $q\ge 2$,\quad $0<t<2/3$, and
\[
t(4-3t)(-K_X^3)\le -4K_X\cdot c_2(X).
\]
\end{parag}

\begin{parag}{$\EEE$ is not semistable, $\rk \FFF=1$, and 
$(\EEE/\FFF)^{**}$ is semistable.}
\label{mcase-1-1}
Then $q\ge 4$,\quad $0<t<1/3$, and
\[
(1-t)(1+3t)(-K_X^3)\le -4K_X\cdot c_2(X).
\]
\end{parag}

\begin{parag}{$\EEE$ is not semistable, $\rk \FFF=1$, and 
$(\EEE/\FFF)^{**}$ is not semistable.}
\label{mcase-1-2}
Then again $q\ge 4$ and $0<t<1/3$. There exists an unstable reflexive sheaf 
$\FFF\varsubsetneq \GGG\varsubsetneq \EEE$. Write
$c_1(\GGG/\FFF)\equiv -p'L$,\ $p'\in \ZZ$ and put $u:=p'/q$, so that
$c_1(\GGG/\FFF)\equiv uK_X$. Then $t<u<1-t-u$ and 
\[
\Bigl(tu+(t+u)(1-t-u)\Bigr)(-K_X^3)\le -K_X\cdot c_2(X),
\]
\end{parag}
\end{mparag}

\begin{mtparag}{Corollary.}
\label{cor-mcase-index}
If $q\QQ(X)=1$, then $\EEE$ is semistable.
If $q\QQ(X)\le 3$, then either $\EEE$ is semistable or 
we are in case \xref{mcase-2}.
\end{mtparag}

\section{Two birational constructions}
\label{sect-3}
\begin{mparag}{}
\label{assumption-Sarkisov-link-0}
Let $X$ be a $\QQ$-Fano threefold. 
Throughout this paper we assume that the linear 
system $|-K_X|$ is non-empty, has no fixed components,
and is not composed of a pencil.
Then a general member $H\in |-K_X|$ is irreducible.
By \eqref{eq--K-exact-2} and \ref{prop-A-lin-syst}
this holds automatically when $-K_X^3\ge 8$. 
Let $q:=q\QQ(X)$ and 
$L$ be the ample Weil divisor that generates the group 
$\Cl X/\mathrm{Tors}$.
Thus we have $-K_X\equiv q L$.
Put $\HHH:=|-K_X|$. Let $H\in \HHH$ be a general member.
\end{mparag}

\begin{mparag}{}
\label{assumption-Sarkisov-link}
Assume there is a diagram \textup(Sarkisov link of type I or II\textup)
\begin{equation}
\label{eq-main-diag-1}
\xymatrix{
\tilde X\ar[d]_{g}\ar@{-->}[r]^{\chi}&Y\ar[d]^{f}
\\
X&Z}
\end{equation}
where $\tilde X$ and $Y$ have only $\QQ$-factorial terminal
singularities, $\rho(\tilde X)=\rho (Y)=2$, 
$g$ is a Mori extremal divisorial contraction,
$\tilde X \dashrightarrow Y$ is a sequence of 
log flips, and $f$ is a Mori extremal contraction
(either divisorial or fibre type). 
Thus one of the following holds:
a) $\dim Z=1$ and $f$ is a $\QQ$-del Pezzo fibration,
b) $\dim Z=2$ and $f$ is a $\QQ$-conic bundle, or
c) $\dim Z=3$, $f$ is a divisorial contraction, and 
$Z$ is a $\QQ$-Fano.
Let $E$ be the 
$g$-exceptional divisor. 
We assume that 
the composition $f{\comp} \chi {\comp} g^{-1}$
is not an isomorphism. For a divisor $D$ on $X$,
everywhere below $\tilde D$ and $D_Y$ denote 
strict birational transforms of $D$ on $\tilde X$ and $Y$, respectively.
We also assume that the discrepancy 
$\alpha:=a(E,X,\HHH)$ is non-positive, i.e.,
\begin{equation}
\label{eq-crepant-formula}
0\sim f^*(K_X+\HHH) =K_{\tilde X}+\tilde \HHH+\alpha E,\quad \alpha \in \ZZ,\quad \alpha \ge 0.
\end{equation} 
By the above we have
\begin{equation}
\label{eq-crepant-formula-corollary}
\dim |-K_{\tilde X}|\ge \dim \tilde \HHH= \dim |-K_X|.
\end{equation}
\end{mparag}

\begin{mparag}{}
\label{eq-beliberda-beta}
Similarly, 
\[
0\qq g^*(K_X+qL)\qq K_{\tilde X}+q\tilde L+ \beta E.
\]
Therefore,
\begin{equation}
\label{eq-beliberda-3}
K_{Y}+qL_{Y}+\beta E_{Y}\qq 0.
\end{equation}
If $q\QQ(X)=qW(X)$, then $K_X+qL\sim 0$ and $\beta$ is an integer
$\ge \alpha$.

Let $F=f^{-1}(\mathrm{pt})$ be a general fibre.
Recall that $F$ is either $\PP^1$ or a smooth del Pezzo surface.
Restricting \eqref{eq-beliberda-3} to $F$
we get 
\begin{equation}
\label{eq-beliberda-4}
K_{F}+qL_{Y}|_{F}+\beta E_{Y}|_{F}\sim 0.
\end{equation}
Here $-K_F$, $L_{Y}|_{F}$, and $E_{Y}|_{F}$ are proportional
nef Cartier divisors. Moreover, $-K_F$ and $E_{Y}|_{F}$
are ample.
\end{mparag}

\begin{mparag}{}
We will use construction \eqref{eq-main-diag-1} in the 
following two situationa:
\begin{parag}{} 
\label{construction-1a}
(see \cite{Takagi-2002-I-II}, \cite{Takagi-2006}).
Let $P\in X$ be a singularity of index $r$.
Take $g$ to be a divisorial blowup of $P$ such that
the discrepancy of the exceptional divisor $E$ is equal to $1/r$.
Assume that the divisor $-K_{\tilde X}$ is nef, big
and the linear system $|-nK_{\tilde X}|$ does not contract any divisors.
Then the transformation in 
\eqref{eq-main-diag-1} is so-called ``two rays game''. 
If $-K_{\tilde X}$ is ample, then $f {\comp} \chi$
is a composition of steps of the $K$-MMP.
Otherwise, $f {\comp} \chi$ is a composition of a single flop 
followed by steps of the $K$-MMP.
It is easy to see 
also that $f {\comp} \chi$ is an $-E$-MMP. 
\end{parag}

\begin{parag}{} 
\label{construction-1b}
(see \cite{Alexeev-1994ge}).
The pair $(X,\HHH)$ is not canonical.
Let $c$ be the canonical threeshold of $(X,\HHH)$.
Then $0<c<1$.
Take $g$ to be an extremal divisorial $K_X+c\HHH$-crepant
blowup.
In this situation, $\alpha>0$ and $f {\comp} \chi$ is an $K+c\HHH$-MMP.
In particular, $f$ is an extremal $K_X+c\HHH$-negative contraction.
The conditions of \ref{assumption-Sarkisov-link} are satisfied by
\cite{Alexeev-1994ge}.
\end{parag}
\end{mparag}

\begin{mparag}{Properties of construction \ref{assumption-Sarkisov-link}.}

\begin{tparag}{Claim.}
\label{claim-E-mom-except}
$E_Y$ is not contracted by $f$.
\end{tparag}
\begin{proof}
Assume the converse, i.e., $\dim f(E_Y)<\min (2, \dim Z)$.
If $f$ is birational, this implies that the map
$f{\comp} {\chi} {\comp} g^{-1}\colon X \dashrightarrow Z$ is an
isomorphism in codimension one.
Since both $X$ and $Z$ are Fano threefolds, this implies that
$f{\comp} \chi {\comp} g^{-1}$ is in fact an
isomorphism. This contradicts our assumptions.
If $\dim Z\le 2$, then $E_Y$ is a pull-back of an ample Weil 
divisor on $Z$. But then $nE_Y$ is movable for some 
$n>0$. Again we derive a contradiction.
\end{proof}

\begin{tparag}{Claim.}
\label{lemma--K-nef}
For some $n,\, m>0$ there is a decomposition 
$-nK_{\tilde X}\sim m\tilde \HHH+ M$, where $|M|$ is a base point free
linear system. In particular, $|-nK_{\tilde X}|$ has no fixed components.
\end{tparag}
\begin{proof}
By \eqref{eq-crepant-formula}, for some $0<c\le 1$, we have 
$K_{\tilde X}+c\tilde \HHH= g^*(K_X+cH)$.
Hence we can take $n,\, m>0$ so that 
$|-nK_{\tilde X}-m\tilde \HHH|$ is base point free.
\end{proof}

\begin{tparag}{Lemma (\cite{Mori-Prokhorov-2006}).}
If $f$ is a $\QQ$-conic bundle, then $Z$ is a del Pezzo surface with at worst Du Val 
singularities of type $A_n$ and $\rho(Z)=1$.
Moreover, there is a natural embedding $f^*\colon \Cl Z\to \Cl Y$.
\end{tparag}
\begin{proof}
The assertion about the base is an immediate
consequence of the main result of \cite{Mori-Prokhorov-2006}
and the fact that $Z$ is uniruled.
The last statement is obvious because both $Y$ and $Z$ 
have only isolated singularities and $\Pic (Y/Z)\simeq \ZZ$.
\end{proof}

\begin{tparag}{Remark.}
\label{rem-class-delPezzo}
\textup{(i)}
In the above notation the generic fibre of $f$ is 
a smooth rational curve. The locus 
$\Lambda:=\{z \in Z \mid \text{$f$ is smooth over $z$}\}$
is a closed subset of codimension $\ge 1$ in $Z$.
The union of one-dimensional components of $\Lambda$
is called the \emph{discriminant curve}.

\textup{(ii)}
The classification of del Pezzo surfaces $Z$ with Du Val 
singularities and $\rho(Z)=1$ is well-known.
In particular, we always have $K_Z^2\le 9$ and $K_Z^2\neq 7$. 
Moreover,
\begin{enumerate}
\item 
if $K_Z^2=9$, then $Z\simeq \PP^2$;
\item 
if $K_Z^2=8$, then $Z\simeq \PP(1,1,2)$;
\item 
if $K_Z^2\le 6$, then on $Z$ there is a rational curve 
$C$ such that $-K_Z\cdot C=1$.
\end{enumerate}
\end{tparag}

\begin{tparag}{Lemma.}
\label{lemma-constr-diag-descr}
Notation and assumptions as in \xref{assumption-Sarkisov-link}.
Assume additionally that $\dim |L|>0$, $q\QQ(X)\ge 4$ and $f$ is not birational.
Then $L_Y=f^*\Xi$ for some \textup(integral\textup) Weil divisor on $Z$.
Moreover, $\dim |\Xi|=\dim |L|$ and the class of $\Xi$ generates the group 
$\Cl Z/\mathrm{Tors}$.
\end{tparag}
\begin{proof}
Since $q\QQ(X)\ge 4$, relation \eqref{eq-beliberda-4} implies $L_{Y}|_{F}=0$.
Since $f$ is a Mori contraction and $Z$ is normal, $L_Y=f^*\Xi$, where $\Xi:=f(L_Y)$.
The rest is obvious.
\end{proof}

\begin{tparag}{Lemma.}
\label{lemma-constr-diag-dimZ2-not-canon}
Assume that $(X,|-K_X|)$ is not canonical
and we are applying construction \xref{assumption-Sarkisov-link}. 
Further, assume that $\dim Z= 2$ and $\alpha>0$.
Then one of the following holds:
\begin{enumerate}
\item 
$\HHH_Y$ is $f$-ample. Then the discriminant curve of $f$ 
is empty.
\item 
$\HHH_Y$ is not $f$-ample. Then $q\QQ(X)\ge 7$.
Moreover, the equality holds only if $Z\simeq \PP^2$
and $\dim |-K_X|=35$.
\end{enumerate}
\end{tparag}

\begin{proof}
First we assume that $\HHH_Y$ is $f$-ample. 
By \eqref{eq-crepant-formula} and Claim \ref{claim-E-mom-except}
$E_Y$ and general elements of $\HHH_Y$ are sections of $f$.
Hence $f$ is smooth outside of 
a finite number of degenerate fibres.

Now we assume that $\HHH_Y$ is not $f$-ample. 
Then $\HHH_Y=f^* \MMM$, where $\MMM$ is a linear system
without fixed components. Let $\Xi$ be an ample Weil
divisor that generates $\Cl Z/\mathrm{Tors}$.
We can write $\MMM\qq a \Xi$
and $-K_Z\qq q' \Xi$, where $q':=q\QQ(Z)$, $a\in \ZZ$.
Clearly, $q\QQ(X)\ge a$.

By our assumption and by Reid's Riemann-Roch formula \cite[(9.1)]{Reid-YPG1987},
\[
30\le \dim \MMM\le \frac 12 \MMM \cdot (\MMM-K_Z)+\sum c_P(\MMM)
\le \frac { a( a+q')}{2q'^2} K_Z^2.
\]
Assume that $a\le 7$. 
If $K_Z^2\le 6$, then $q'=K_Z^2$ by Remark \ref{rem-class-delPezzo}.
So, 
$60q'\le a( a+q')\le 49+7q'$, a contradiction.
If $K_Z^2=8$, then $q'=4$, so $120\le a( a+4)\le 77$. Again 
we have 
a contradiction. Finally, let $K_Z^2=9$, i.e.,
$Z\simeq\PP^2$. Then $q'=3$, so $60\le a( a+3)\le 70$.
This inequality has only one solution: $a=7$.
But then $q\QQ(X)\le 7$. If $q\QQ(X)= 7$,
then $a=7$, $\MMM=|\OOO_{\PP^2}(7)|$,
and $\dim \MMM=35$.
\end{proof}

\begin{tparag}{Lemma.}
\label{lemma-constr-diag-dimZ2-q1}
Notation and assumptions as in \xref{assumption-Sarkisov-link}.
Assume additionally that $q\QQ(X)=1$, $Z$ is a surface, 
and the discriminant curve of 
$f$ is empty.
Then $\dim |-K_X|<30$.
\end{tparag}
\begin{proof}
Suppose $\dim |-K_X|\ge 30$.
Let $\Gamma\subset Z$ is a smooth curve contained into the 
smooth locus of $Z$. Then 
$G:=f^{-1}(\Gamma)$ is a smooth ruled surface over $\Gamma$.
We claim that $\dim |-K_Y-G|\le 0$.
Indeed, otherwise $-K_Y\sim G+B$, where $B$ 
is an integral effective divisor, $\dim |B| \ge 1$. 
Since $q\QQ(X)=1$, this gives a contradiction.

Now from \eqref{eq-crepant-formula-corollary} and from the exact sequence 
\[
0\longrightarrow \OOO_Y(-K_Y-G)
\longrightarrow \OOO_Y(-K_Y)
\longrightarrow \OOO_G(-K_Y)
\longrightarrow 0
\]
we get $h^0(\OOO_G(-K_Y))\ge h^0(\OOO_Y(-K_Y))-1\ge 30$.
It is easy to see that 
\[
(-K_Y|_G)^2=(-K_G+G|_G)^2=K_G^2-2K_G\cdot G|_G=8-8p_a(\Gamma)+4\Gamma^2.
\]
By Claim \ref{lemma--K-nef} the linear system $|-nK_Y|$ has 
no fixed components. Therefore we can take $\Gamma$ so that
$|-nK_Y|_G|$ has at worst isolated base points 
(in particular, it is nef). 
Moreover, $|-nK_Y|_G|$ is base point free for sufficiently
large $n$. If $-K_Y|_G$ is ample, 
it is well-known that $h^0(\OOO_G(-K_Y))\le (-K_Y|_G)^2+2$
(see, e.g., \cite{Fujita-1975}). 
If $-K_Y|_G$ is not ample, we 
obtain the above inequality by 
applying the same arguments to 
$\bar G$, where $\bar G$ is the image of $G$ under the birational
contraction given by $|-nK_Y|_G|$.
In both cases we have 
\[
8-8p_a(\Gamma)+4\Gamma^2=(-K_Y|_G)^2\ge h^0(\OOO_G(-K_Y)) -2\ge 28.
\]
This gives us 
\[
\Gamma^2\ge 2p_a(\Gamma)+5=K_Z\cdot \Gamma+\Gamma^2+7,
\qquad -K_Z\cdot \Gamma\ge 7.
\]
If $K_Z^2<8$, then we can take $\Gamma$ to be a general member
of $-K_Z$ and derive
a contradiction. If $K_Z^2=8$ or $9$, then we can take 
$\Gamma \in |-\frac 12 K_Z|$, or $|-\frac 13 K_Z|$, respectively.
\end{proof}

\begin{tparag}{Lemma.}
\label{cor-dim1-index3}
If $\dim Z=1$ and $\dim |-K_X|\ge 30$, then $q\QQ(X)\ge 3$.
\end{tparag}

\begin{proof}
Let $F_1,F_2, F_3$ be general fibres.
Then from the exact sequence
\[
0 \longrightarrow
\OOO_{Y}\left(-K_Y-\sum F_i\right)\longrightarrow\OOO_Y(-K_Y)
\longrightarrow \bigoplus \OOO_{F_i}(-K_{F_i}) \longrightarrow 0
\]
we obtain
\[
h^0(-K_Y-\sum F_i)\ge h^0(-K_Y)- \sum
h^0(-K_{F_i}).
\]
Since $F_i$ are smooth del Pezzo surfaces,
$h^0(-K_{F_i})=K_F^2+1\le 10$.
Hence, $h^0(-K_Y-\sum F_i)>0$ by \eqref{eq--K-exact-2} and
we have a decomposition $-K_Y\sim \sum F_i+G$, where $G$ is effective. 
Since $F_i$ is movable, this gives us that $q\QQ(X)\ge 3$.
\end{proof}
\end{mparag}

\begin{mparag}{Case: $(X,|-K_X|)$ is canonical.}
\begin{parag}{}
\label{assumpt-ge}
Consider the case when $(X,|-K_X|=\HHH)$ is canonical.
According to \cite{Alexeev-1994ge} there is the following
diagram
\[
\xymatrix{&\tilde X\ar[dr]^f\ar[dl]_g\ar[r]&\bar X\ar[d]&
\\
X\ar@{-->}[rr]&&Y\subset \PP^n}
\]
where $g\colon (\tilde X,\tilde \HHH)\to (X,\HHH)$ is a 
terminal modification of $(X,\HHH)$, $n:=\dim |-K_X|$,
the morphism $f$ is given by 
the (base point free) linear system $\tilde \HHH$,
$\dim Y=2$ or $3$, and 
$\tilde X\to \bar X\to Y$ is the Stein factorization.
We have
\[
K_{\tilde X}+\tilde \HHH=g^*(K+\HHH)\sim 0.
\]
Since $(\tilde X,\tilde \HHH)$ is terminal,
a general member $\tilde H\in \tilde \HHH$ is a smooth K3 surface.
From the exact sequence
\[
0\longrightarrow
\OOO_{\tilde X}\longrightarrow
\OOO_{\tilde X}(-K_{\tilde X})\longrightarrow
\OOO_{\tilde H}(-K_{\tilde X})\longrightarrow
0
\]
one can see that the restriction 
$f|_{\tilde H}$ is given by a complete linear system. 
\end{parag}

\begin{tparag}{Lemma.}
\label{lemma-canonical-2}
Let $X$ be a $\QQ$-Fano threefold.
Assume that $(X,|-K_X|=\HHH)$ is canonical
and the image of the map given by $|-K_X|$
is a surface. If $\dim |-K_X|\ge 6$, then
$2q\QQ(X)\ge \dim |-K_X|-1$.
\end{tparag}

\begin{proof}
We use notation of \ref{assumpt-ge}.
By our assumption 
$f(\tilde H)$ is a curve. Thus $|-K_{\tilde X}|_{\tilde H}|$
is a base point free elliptic pencil 
on $\tilde H$ and $f(\tilde H)\subset \PP^n$
is a rational normal curve of degree $n-1$.
Hence $Y\subset \PP^{n}$ is a surface of degree $n-1$.
Let $M$ be a hyperplane section of $Y$.
It is well-known 
that 
in this situation one of the following halds
(recall that $n\ge 6$): 
\begin{enumerate}
\item 
$Y$ is a rational scroll, $Y\simeq\FF_e$, $M\sim \Sigma+al$, 
where $\Sigma$ and $l$ are the minimal section and 
a fibre of $\FF_e$, respectively, and $a$ is an integer such that
$a\ge e+1$, $n-1=2a-e$.

\item 
$Y$ is a cone over a rational normal curve of degree $n-1$,
$M\sim (n-1)l$, where $l$ is a generator of the cone.
\end{enumerate}
In case (i), $\tilde \HHH \sim f^*\Sigma+af^*l$. Here
$|f^*l|$ is a linear system 
without fixed components and $f^*\Sigma$ is an effective divisor.
So, $2q\QQ(X)\ge 2a\ge n-1$.
In case (ii) we have $\tilde \HHH \sim f^*(n-1)l$.
Let $o\in Y$ be the vertex of the cone and let
$G$ be the closure of $f^*l$ over $Y\setminus \{o\}$.
Then $G$ is an integral Weil divisor and $\tilde H\qq (n-1)G+ T$,
where $T$ is effective.
Clearly, $g$ does not contract any component of $G$.
This implies $q\QQ(X)\ge n-1$.
\end{proof}

Now assume that $\dim Y=3$.
\begin{tparag}{Lemma (cf. {\cite[Corollary 1.8]{Prokhorov-2005a}}).}
\label{lemma-canonical-3}
Let $X$ be a $\QQ$-Fano threefold.
Assume that $(X,|-K_X|=\HHH)$ is canonical
and the image of the map given by $|-K_X|$
is three-dimensional. 
Then $\dim |-K_X|\le 37$. 
If moreover $q\QQ(X)=1$, then $\dim |-K_X|\le 13$.
\end{tparag}

\begin{proof}
By the construction, 
$\bar Y$ is a Fano threefold with canonical Gorenstein 
singularities 
and $\bar Y\to Y\subset \PP^N$ is the anticanonical map
(see \cite{Alexeev-1994ge}).
We have $\dim |-K_X|\le \dim |-K_{\bar Y}|\le 38$ 
by the main result of \cite{Prokhorov-2005a}.
Moreover, if $\dim |-K_X|=38$, then 
$\bar Y$ is isomorphic either
$\PP(3,1,1,1)$ or $\PP(6,4,1,1)$.
In particular, $\bar Y$ is a toric variety.
Since $\tilde X$ is a terminal modification
of $\bar Y$, it is also toric and so is $X$.
By Lemma \ref{lemma-toric} below $\dim |-K_X|\le \dim |-K_{\bar Y}|\le 33$,
a contradiction.
If $q\QQ(X)=1$, 
then $-K_{\bar Y}$ cannot be decomposed into a 
sum of two movable divisors. 
According to \cite{Mukai-2002}, $\dim |-K_X|\le \dim |-K_{\bar Y}|\le 13$.
\end{proof}

\begin{tparag}{Lemma.}
\label{lemma-toric}
Let $X$ be a toric $\QQ$-Fano threefold. 
If $X\not\simeq\PP^3$, then $-K_X^3\le 125/2$ and $\dim |-K_X|\le 33$. 
\end{tparag}

\begin{proof}[Sketch of the proof]
By considering cyclic covering tricks
(cf. Proof of Proposition \ref{proposition-Cl-2})
we reduce the question to the case $\Cl X\simeq \ZZ$.
For toric varieties this preserves the property 
$\rho=1$. Then $X$ is a weighted projective space.
Using the fact that $X$ has only terminal singularities
we get the following cases:
$\PP(1, 1, 1, 2)$, $\PP(1, 1, 2, 3)$, $\PP(1, 2, 3, 5)$, $\PP(1, 3, 4, 5)$,
$\PP(2, 3, 5, 7)$, $\PP(3, 4, 5, 7)$.
The lemma follows.
\end{proof}
\end{mparag}

\section{Case $q\QQ(X)\le 3$}
\label{sect-4}
In this section we consider the case $q:=q\QQ(X)\le 3$.

\begin{mtparag}{Proposition.}
\label{cases-I}
Let $X$ be a $\QQ$-Fano threefold.
Assume that $X$ is not Gorenstein, $q:=q\QQ(X)\le 3$ and $-K_X^3\ge 125/2$. 
Then we have one of the following cases:
\begin{parag}{}
\label{cases-I-2}
$q=1$, $\B=(2)$, $-K_X^3=2g-3/2$, $\dim |-K_X|=g+1$, $32\le g\le 35$;
\end{parag}

\begin{parag}{}
\label{cases-I-2-2}
$q=1$, $\B=(2,2)$, $-K_X^3= 63$, $\dim |-K_X|=33$;
\end{parag}

\begin{parag}{}
\label{cases-I-3-1}
$q=1$, $\B=(3)$, $-K_X^3=188/3$, $\dim |-K_X|=33$;
\end{parag}

\begin{parag}{}
\label{cases-I-3-2}
$q=2$, $\B=(3)$, $L^3=25/3$, $\dim |L|=9$, $\dim |-K_X|=35$.
\end{parag}
\end{mtparag}

\begin{mtparag}{Lemma.}
\label{lemma-qQ3}
In notation of Proposition \xref{cases-I} we have 
$-K_X\cdot c_2(X)\ge 125/8$ and
$\sum_{P\in \B} (r_P-1/r_P)\le 67/8$. In particular,
$\sum r_P\le 10$.
\end{mtparag}
\begin{proof}
By Corollary \ref {cor-mcase-index}
we have cases \ref{mcase-semi} or \ref{mcase-2}.
Hence,
\[
-K_X\cdot c_2(X)\ge 
\begin{cases}
\displaystyle{\frac13 (-K_X^3)\ge \frac{125}6}, 
\\[10pt]
\displaystyle{
\frac 14 t(4-3t) (-K_X)^3\ge \frac 1{4q} \left(4-\frac 3q\right) 
\frac{125}2\ge \frac{125}8. 
}
\end{cases}
\]
(In the second line we used that $t\ge 1/q\ge 1/3$
and the function $t(4-3t)$ is increasing for
$t\le 2/3$).
In both cases we have $-K_X\cdot c_2(X)\ge 125/8$.
Thus,
\[
\sum_{P\in \B} \left( r_P-\frac1{r_P}\right) \le 24-\frac{125}{8}=\frac{67}{8}.
\]
Hence $\B$ contains at most $5$ points and 
$\sum r_P\le \down{\frac{67}{8}+5\cdot \frac 12}\le 10$. 
\end{proof}

\begin{mtparag}{Proposition.}
\label{proposition-Cl-1}
In notation of Proposition \xref{cases-I} we have 
$\Cl X\simeq \ZZ$.
\end{mtparag}
\begin{proof}
Let $T$ be an $s$-torsion element in the Weil divisor
class group. By Riemann-Roch \eqref{eq-RR}, Kawamata-Viehweg vanishing theorem and
Serre duality we have
\[
\begin{array}{lll}
0&=\chi(T)&=1+\sum _P c_P(T),
\\[10pt]
0&=\chi(K_X+T)&=1+\frac1{12}K_X\cdot c_2(X)+\sum_{P\in \B} c_P(K_X+T).
\end{array}
\]
Subtracting we get
\[
0=-\frac1{12}K_X\cdot c_2(X)+\sum_{P\in \B} (c_P(T)-c_P(K_X+T)).
\]
Take $i_{T,P}$ so that $T\sim i_{T,P}K_X$ near $P\in \B$.
Then $si_{T,P}\equiv 0\mod r_P$ and
\[
0=-\frac1{12}K_X\cdot c_2(X)+\frac {1}{12}
\sum_{P\in \B} \left( r_P-\frac{1}{r_P}\right) -
\sum_{P\in \B}
\frac{\ov{b_Pi_{T,P}}\ \left(r_P-\ov{b_Pi_{T,P}}\right)}{2r_P}.
\]
Therefore,
\[
2=\sum_{P\in \B}
\frac{\ov{b_Pi_{T,P}}\ \left(r_P-\ov{b_Pi_{T,P}}\right)}{2r_P}.
\]
If $i_{T,P}\not \equiv 0\mod r_P$, we have 
\[
\frac{\ov{b_Pi_{T,P}}\ \left(r_P-\ov{b_Pi_{T,P}}\right)}{2r_P}\le 
\frac {r_P}8.
\]
Combining the last two relations we get
\[
\sum_{P\in \B'} r_P\ge 16,
\]
where the sum runs over all $P\in \B$
such that $i_{T,P}\not \equiv 0\mod r_P$. This contradicts Lemma \ref{lemma-qQ3}.
\end{proof}

\begin{proof}[Proof of Proposition \xref{cases-I}]
By Proposition \xref{proposition-Cl-1}
$q=q\QQ(X)=qW(X)$. So,
$\gcd(q,r_P)=1$ for all $P\in \B$. 

\begin{mparag}{Case $q=3$.}
We will show that this case does not occur.
By \eqref{eq-comput-L3} we have 
\begin{equation}
\label{eq-vspom-q=3}
-K_X^3=q^3L^3=162 +\frac 92 K_X\cdot c_2(X)+162 \sum_{P\in \B} c_P(-L). 
\end{equation}
By Lemma \ref{lemma-qQ3} $-K_X\cdot c_2(X)\ge 125/8$ and
$-K_X^3\ge 125/2$ by our assumptions. Combining this we obtain 
$\sum c_P(-L)\ge -467/2592$.

Again by Lemma \ref{lemma-qQ3} we have 
$\sum (r_P-1/r_P)\le 67/8$.
Assume that $r_P=2$ for all $P\in \B$. 
Note that $c_P(L)=-1/8$ (because $-K_X\sim L$ near each $P$). 
Hence $\B=(2)$. 
Then 
$-K_X\cdot c_2(X)=45/2$.
By \eqref{eq-vspom-q=3} we have $-K_X^3=81/2<125/2$,
a contradiction. 

Thus we assume that 
at least one on the $r_P$'s is $\ge 3$. 
Recall that 
$\sum r_P\le 10$, 
$\sum (r_P-1/r_P)\le 67/8$
and $3 \nmid r_P$. 
This gives us the following possibilities
for $\B$:
\[ 
(4), (5), (7), (8),
(2, 4), (2, 5), (2, 7),
(2, 2, 4), (2, 2, 5),
(4, 4), (2, 2, 2, 4).
\]
Take $0\le i_P<r_P$ so that $3i_P \equiv -1 \mod r_P$.
Easy computations give us
\[
\begin{array}{|c||c|c|c|c|c|}
\hline
&&&&&
\\
r_P&2&4&5&7&8
\\[8pt]
\hline
&&&&&
\\
i_P&1&1&3&2&5
\\[8pt]
\hline
&&&&&
\\
c_P&-1/8&-5/16&-1/5&-2/7,-3/7, -5/7&-5/32
\\[9pt]
\hline
\end{array}
\]
In all cases except for $\B=(8)$ 
we get a contradiction with $\sum c_P(-L)\ge -467/2592$.
Consider the case $\B=(8)$.
Then by \eqref{eq-vspom-q=3} we have
\begin{equation*}
-K_X^3=162 -\frac 92 \cdot \frac{129}8-162 \frac5{32}=
\frac{513}8. 
\end{equation*}
Then by \eqref{eq--K-exact}
\[
\dim |-K_X|=2+\frac{513}{16}-\frac{b_P(8-b_P)}{16}=
34+\frac{1-b_P(8-b_P)}{16}.
\]
This number cannot be an integer,
a contradiction. 
\end{mparag}

\begin{mparag}{Case $q=1$.}
By \ref{mcase-semi} we have 
\[
\sum_{P\in \B} \left(r_P-\frac1 {r_P}\right)=24+K_X\cdot c_2(X)\le 24+\frac 12 K_X^3
\le 24 -\frac{125}{6}=\frac{19}{6}.
\]
This gives the following possibilities:
$\B=(2)$, $(3)$, or $(2,2)$.

If $\B=(2,2)$, then $-K_X\cdot c_2(X)=21$ and $-K_X^3\le 63$. 
On the other hand, 
$-K_X^3\in \frac12 \ZZ$ (see \cite[Lemma 1.2]{Suzuki-2004}). 
Hence $-K_X^3= 63$ or $125/2$.
Further, by \eqref{eq--K-exact}
\[
\dim |-K_X|=-\frac 12 K_X^3+\frac {3}{2}. 
\]
Since this number should be an integer, the only possibility is $-K_X^3= 63$
and $\dim |-K_X|=33$.

If $\B=(2)$, then $-K_X\cdot c_2(X)=45/2$ and by \eqref{eq--K-exact}
\[
\dim |-K_X|
= -\frac 12 K_X^3 +\frac{7}{4}. 
\]
Put $g:=\dim |-K_X|-1$. Then $-K_X^3=2g - 3/2$.
We have
\[
125/2 \le -K_X^3=2g - 3/2\le 74-\frac 92.
\]
Hence $32\le g\le 35$ and $-K_X^3\in\{125/2, 129/2, 133/2, 137/2\}$.

Assume that $\B=(3)$. Then $-K_X\cdot c_2(X)=64/3$
and $-K_X^3\le 64$. 
As above,
\[
\dim |-K_X|=
-\frac 12 K_X^3 +\frac53. 
\]
We get only one possibility:
$-K_X^3=188/3$ and $\dim |-K_X|= 33$.
\end{mparag}

\begin{mparag}{Case $q=2$.}
If $\EEE$ is semistable, then as above by \ref{mcase-semi} $\B=(3)$.
Otherwise we are in case \ref{mcase-2} and as in the proof of Lemma 
\ref{lemma-qQ3} we have 
\[
\sum_{P\in \B} \left(r_P-\frac1 {r_P}\right)=24+K_X\cdot c_2(X)\le 24
+\frac 5{16}K_X^3\le \frac{143}{32}.
\]
Since $\gcd(r_P,q)=1$, again we get the same possibility $\B=(3)$.

Then $-K_X\cdot c_2(X)=64/3$ and $L\cdot c_2(X)=32/3$. Hence
\[
5/4(-K_X^3)\le t(4-3t)(-K_X^3)\le 4 \cdot 64/3.
\]
Thus $125/2\le -K_X^3\le 1024/15$ and 
$125/16 \le L^3\le 128/15$. Since $3L^3\in \ZZ$ 
(see \cite[Lemma 1.2]{Suzuki-2004}), we have 
$L^3=8$ or $25/3$. 
As above the case 
$L^3=8$ is impossible by \eqref{eq--K-exact}.
Thus $L^3=25/3$.
Then one can easily compute $h^0(L)$
and $h^0(-K_X)$ by \eqref{eq-comput-chi}.
\end{mparag}
\end{proof}

\section{Case $q\QQ(X)\ge 4$}
\label{sect-5}
\begin{mtparag}{Proposition}
\label{cases-II}
Let $X$ be a $\QQ$-Fano threefold.
Assume that $X$ is not Gorenstein, $-K_X^3\ge 125/2$, and
$q:=qW(X)=q\QQ(X)\ge 4$. 
Then we have one of the following cases:

\begin{parag}{}
\label{cases-II-q4-B5}
$q=4$, $\B=(5)$, $-K_X^3= 384/5$, $\dim|L|=3$, $\dim |2L|=10$, $\dim |-K_X|=40$;
\end{parag}

\begin{parag}{}
\label{cases-II-q4-B55}
$q=4$, $\B=(5, 5)$, $-K_X^3= 64$, $\dim|L|=2$, $\dim |2L|=8$, $\dim |-K_X|=33$;
\end{parag}

\begin{parag}{}
\label{cases-II-q5-B2}
$q=5$, $\B=(2)$, $-K_X^3= 125/2$, $\dim|L|=2$, $\dim |2L|=6$, $\dim |-K_X|=33$;
\end{parag}

\begin{parag}{}
\label{cases-II-q5-B26}
$q=5$, $\B=(2, 6)$, $-K_X^3= 250/3$, $\dim|L|=2$, $\dim |2L|=7$, $\dim |-K_X|=43$;
\end{parag}

\begin{parag}{}
\label{cases-II-q5-B7}
$q=5$, $\B=(7)$, $-K_X^3= 500/7$, $\dim|L|=2$, $\dim |2L|=6$, $\dim |-K_X|=37$;
\end{parag}

\begin{parag}{}
\label{cases-II-q5-B2236}
$q=5$, $\B=(2, 2, 3, 6)$, $-K_X^3= 125/2$, $\dim|L|=1$, $\dim |2L|=5$, $\dim |-K_X|=32$;
\end{parag}

\begin{parag}{}
\label{cases-II-q6-B57}
$q=6$, $\B=(5, 7)$, $-K_X^3= 2592/35$, $\dim|L|=1$, $\dim |2L|=4$, $\dim |-K_X|=38$;
\end{parag}

\begin{parag}{}
\label{cases-II-q7-B39}
$q=7$, $\B=(3, 9)$, $-K_X^3= 686/9$, $\dim|L|=1$, $\dim |2L|=3$, $\dim |-K_X|=39$;
\end{parag}

\begin{parag}{}
\label{cases-II-q7-B210}
$q=7$, $\B=(2, 10)$, $-K_X^3= 343/5$, $\dim|L|=1$, $\dim |2L|=3$, $\dim |3L|=6$, $\dim |-K_X|=35$.
\end{parag}
\end{mtparag}

\begin{proof}
Let $L$ be a Weil divisor such that $-K_X\sim q L$.
Since $qW(X)=q\QQ(X)$, the group $\Cl X/\mathrm{Tors}$
is generated by $L$. To get our cases
we run a computer program. 
Below is the description of our algorithm.

1)
By \eqref{eq-RR-O} and Theorem \ref{th-Kawamata-Kc} we have 
$\sum_{P\in \B} (1-1/r_P)\le 24$. Hence there is only 
a finite (but very huge) number of possibilities for the basket $\B$.
In each case we know $-K_X\cdot c_2(X)$ from 
\eqref{eq-RR-O}.
Let $r:=\lcm (\{r_P\})$ be the Gorenstein index of $X$.

2)
By Corollary \ref{cor-first-prop-index} $q\le 4r$
and $\gcd (q,r)=1$.
Hence we have only a finite number of possibilities for 
the index $q$.

3) 
In each case we compute $L^3$ and $-K_X^3=q^3L^3$ by formula
\eqref{eq-comput-L3} and check the condition $-K_X^3\ge 125/2$. 
Here, for $D=-L$, the number $i_P$ is uniquely determined by
conditions $qi_P\equiv b_P\mod r_P$ and $0\le i_P< r_P$.

4) 
Next we check Kawamata's inequalities \ref{Kawamatas-ineq},
i.e., we check that at least one of inequalities
\ref{mcase-semi} -- \ref{mcase-1-2}
holds.
In case \ref{mcase-2} we use the fact that the function
$t(4-3t)$ is increasing for $t<2/3$. 
Since $t\ge 1/q$, we have $\frac1q (4-\frac 3q)\le t(4-3t)$
and 
\[
\frac1q \left(4-\frac 3q\right)\left(-K_X^3\right)\le -4K_X\cdot c_2(X).
\]
Similarly, in cases \ref{mcase-1-1} and \ref{mcase-1-2}
we have, respectively,
\[
\left(1-\frac1q\right)\left(1+\frac3q\right)\left(-K_X^3\right)\le -4K_X\cdot c_2(X),
\]
\[
\frac1q \left(2-\frac3{q}\right)\left(-K_X^3\right)\le -K_X\cdot c_2(X).
\]

5) 
Finally, by the Kawamata-Viehweg vanishing theorem
we have $\chi(tL)=h^0(tL)=0$ for $-q<t<0$.
We check this condition by using 
\eqref{eq-comput-chi}.

At the end we get possibilities \ref{cases-II-q4-B5}--\ref{cases-II-q7-B210}.
\end{proof}

\begin{mtparag}{Corollary (cf. {\cite[Remark 2.14]{Suzuki-2004}}).}
\label{cor-degle}
Let $X$ be a $\QQ$-Fano threefold. If $qW(X)=q\QQ(X)$, then
$-K_X^3\le 250/3$.
\end{mtparag}

Now we show that the condition $qW(X)=q\QQ(X)$
in Proposition \ref{cases-II} is satisfied automatically.

\begin{mtparag}{Proposition.}
\label{proposition-Cl-2}
Let $X$ be a $\QQ$-Fano threefold.
Assume that $q:=q\QQ(X)> 3$
and $-K_X^3> 45$. Then $\Cl X\simeq \ZZ$.
\end{mtparag}

\begin{proof}
Assume that the torsion part of $\Cl X$ is non-trivial for some $X$
satisfying the conditions of Proposition \ref{cases-II}.
Take $X$ so that $q\QQ(X)$ is maximal.
Write $K_X+qL\qq 0$, where $L$ is an (ample) integral Weil divisor.
Since $\Cl X$ is finitely generated and
by cyclic covering trick \cite[(3.6)]{Reid-YPG1987},
there is a finite \'etale in codimension one
cover $\pi\colon X'\to X$ such that $\Cl X'$
torsion free. 
Here $K_{X'}+qL'\sim 0$, where $L':=\pi^*L$.
Note that $X'$ has only terminal singularities.
Hence $X'$ is a Fano threefold with terminal singularities
with $qW(X')\ge q$. 
(It is possible however that $X'$ is not $\QQ$-factorial
and $\rho(X')>1$). Denote $n:=\deg \pi$.
Clearly, $-K_{X'}^3=-nK_X^3\ge -2K_X^3$.
Hence $\dim |-K_{X'}|\ge -K_X^2-2>43$.
Let $\sigma\colon X''\to X'$ be a $\QQ$-factorialization.
(If $X'$ is $\QQ$-factorial, we take $X''=X'$).
Run $K$-MMP on $X''$:\quad $\upsilon\colon X'' \dashrightarrow Y$.
At the end we get a Mori-Fano fibre space $f\colon Y\to Z$.
Let $L'':=\sigma^{-1}(L')$ and $L_Y:=\upsilon_*L''$.
Then $-K_Y\sim qL_Y$. 
If $\dim Z>0$, then for a general fibre $F:=f^{-1}(o)$, $o\in Z$
we have $-K_F\sim q L_Y|_F$. This is impossible if $q>3$. 

In the case $\dim Z=0$, $Y$ is a Fano with $\rho(Y)=1$ and $qW(Y)\ge q$.
By our assumption of maximality of $q=q\QQ(X)$ we have 
$q\QQ(Y)=qW(Y)= q$. 
Hence, $-K_Y^3\le 250/3$ by Corollary \ref{cor-degle}. 
By \eqref{eq--K-exact-1} we have $\dim |-K_Y|\le 43$.
Using \eqref{eq--K-exact-2}
we obtain
\[
43\ge \dim |-K_Y|\ge \dim |-K_{X''}|\ge -\frac12 K_{X''}^3-2
\ge -K_X^3-2. 
\]
Thus $-K_X^3\le 45$,
a contradiction.
\end{proof}

\section{Proof of the main theorem}
\label{sect-6}
\begin{mparag}{}
To construct a Sarkisov link such as in
\eqref{eq-main-diag-1},
we need the following result basically due to 
Ambro and Kawachi. 

\begin{tparag}{Proposition (cf. {\cite[Th. 4.1]{Takagi-2002-I-II}}).}
\label{prop-base-point-H}
Let $X$ be a Fano threefold with terminal singularities,
and let $S$ be an ample Cartier divisor proportional to $-K_X$.
Then the linear system $|S|$ is non-empty and 
a general member of $|S|$ is a reduced irreducible
normal surface whose singularities are at worst
log terminal of type T. Moreover, assume that 
$K_X^2\cdot S>1$ and $qF(X)\ge 1/2$. Then 
a general $S\in |S|$ does not pass through 
non-Gorenstein points \textup(and has at worst Du Val 
singularities\textup). 
\end{tparag}
\begin{proof}
According to \cite{Ambro-1999} the pair $(X,S)$
is plt for a general $S\in |S|$. Then 
singularities of $S$ are of type T by \cite{Kollar-ShB-1988}.
Note that the restriction map 
$H^0(\OOO_X(S))\to H^0(\OOO_S(S))$ is surjective.
Let $P\in \Bs |S|$ be a non-Gorenstein point of $X$.
Then $P\in S$ is a log terminal non-Du Val singularity of type T.

Recall that Kawachi's invariant of a normal 
surface singularity $(S,P)$ is defined as $\delta_P:=-(\Gamma-\Delta)^2$,
where $\Delta$ is the codiscrepancy divisor of $(S,P)$
on the minimal resolution $\hat S\to S$ and 
$\Gamma$ is the fundamental cycle on $\hat S$
(see \cite{Kawachi-Masek-1998}).
If $(S,P)$ is a rational singularity, then 
$\delta_P=\Gamma^2-\Delta^2+4$.
Hence in our case Kawachi's invariant $\delta_P$ 
is integral (because $\Delta^2\in \ZZ$, see \cite{Kollar-ShB-1988}).
On the other hand, $0<\delta_P<2$. Thus $\delta_P=1$.
Now we apply the main result of \cite{Kawachi-Masek-1998}
to the linear system $|S|_S|=|K_S-K_X|_S|$.
It follows that there is a curve $C$ on $S$ passing through $P$
and such that $-K_X\cdot C<1/2$. Since $qF(X)\ge 1/2$, this is impossible.
\end{proof}

\begin{tparag}{Proposition.}
\label{prop-base-point-H-2}
In notation of Proposition \xref{prop-base-point-H}
assume additionally that $(2K_X+S)^2\cdot S\ge 5$
and 
$-(2K_X+S)$ is an ample divisor
which is divisible in $\Cl X/{\mathrm{Tors}}$. Then the linear system 
$|-K_X|$ has only isolated base points.
\end{tparag}
\begin{proof}
Denote the restriction $-K_X|_S$ by $D$.
Since $S$ does not pass through non-Gorenstein points,
$D$ is Cartier.
By the Kawamata-Viehweg vanishing the map
\[
H^0(\OOO_{X}(-K_X))\longrightarrow H^0(\OOO_{S}(D))
\] 
is surjective. Thus it is sufficient to show that 
the linear system $|D|$ is base point free.
By the adjunction formula $D=K_S-(2K_X+S)|_S$.
Let $\mu \colon \hat S\to S$ be the minimal resolution.
Since $S$ has at worst Du Val singularities,
$K_{\hat S}=\mu^*K_S$. Thus we can write 
$\mu^*D=K_{\hat S}+M$, where $M=\mu^*(-(2K_X+S)|_S)$ is nef.
It is easy to see that $M^2=(2K_X+S)^2\cdot S\ge 5$
by our assumption.
Suppose that the linear system $|\mu^*D|=|K_{\hat S}+M|$
has a base point $P$. 
By the main theorem of \cite{Reider-1988} 
there is an effective divisor $E$ on $\hat S$ passing through $P$
such that either $M\cdot E=0$, $E^2=-1$ or $M\cdot E=1$, $E^2=0$.
In the former case $E$ is contracted my $\mu$ and we get a 
contradiction by the genus formula. In the latter case we have 
$-(2K_X+S)\cdot \mu(E)=1$. This is impossible because 
$-(2K_X+S)$ is divisible in $\Cl X/{\mathrm{Tors}}$ and
$\mu(E)$ is contained in the Gorenstein locus of $X$.
\end{proof}

Since $qF(X)=q/r$, we have the following

\begin{tparag}{Corollary.}
\label{cor--K-isol-base-points}
Let $X$ be a $\QQ$-Fano threefold,
let $q:=q\QQ(X)$, and let $r$ be the Gorenstein index of $X$. 
Assume that $-K_X^3>q/r=qF(X)$, $2q-r\ge 2$, and
$(-K_X^3)(2q-r)^2r\ge 5q^3$. 
Then the linear system 
$|-K_X|$ has only isolated base points.
\end{tparag}

\begin{proof}
Let $L$ be the Weil divisor such that 
$-K_X\qq qL$.
Take $S=rL$ and apply Proposition \ref{prop-base-point-H-2}.
\end{proof}
\end{mparag}

Now we are in position to prove Theorem \ref{main-theorem}.

\begin{mparag}{Main assumption.}
\label{main-assumpt}
Let $X$ be a $\QQ$-Fano threefold. We assume that 
$-K_X^3\ge 125/2$.
Then
 $X$ is such as in Propositions \ref{cases-I}
or \ref{cases-II}.
In particular, 
$\dim |-K_X|\ge 32$.
By Propositions \ref{proposition-Cl-1} and \ref{proposition-Cl-2}
we also have $\Cl X\simeq \ZZ$.
We divide cases of \ref{cases-I}
or \ref{cases-II} in four groups and treat these 
groups separately (see \ref{subsect-case-last-real},
\ref{case-last-mnogo} \ref{case-last-5.1.7}, \ref{case-last-index1}). 
\end{mparag}

\begin{tparag}{Proposition.}
\label{prop-diag-0}
Notation and assumptions as in \xref{main-assumpt}.
If there exists a Sarkisov link \eqref{eq-main-diag-1} 
with birational $f$, then $-K_Z^3\ge 125/2$ except possibly for 
the following case
\begin{itemize}
\item 
$\dim |-K_Z|= \dim |-K_X|= 32$.
\end{itemize}
\end{tparag}

\begin{proof}
Assume the converse.
Then $Z$ is a $\QQ$-Fano with $\dim |-K_Z|\ge \dim |-K_X|\ge 32$
and $-K_Z^3< 125/2$.
By \eqref{eq--K-exact-1} 
\begin{equation}
\label{eq-KZ-K3}
\dim |-K_Z|+\frac 12 \sum_{P\in \B_Z} \left(1- \frac1 {r_P}\right)\le
-\frac 12 K_Z^3 +
2 < \frac{133}4.
\end{equation}
Therefore, $\dim |-K_Z|=32$ or $33$.
Moreover, if $\dim |-K_Z|=33$, then we have $r_P=1$ for all $P\in \B_Z$,
i.e., $Z$ is Gorenstein (and factorial). 
In particular, $q\QQ(Z)=qF(Z)=qW(Z)$ and $q\QQ(Z)^3$
divides $-K_Z^3$.
By Riemann-Roch, $-K_Z^3=62$.
Therefore, $q\QQ(Z)=1$.
But then $-K_Z$ cannot be decomposed into a sum of
movable divisors. 
We derive a contradiction by \cite{Mukai-2002}. 
\end{proof}

\begin{mparag}{Case \xref{cases-II-q5-B2}}
\label{subsect-case-last-real}
\begin{tparag}{Proposition (see {\cite{Sano-1996}}).}
\label{prop-case-last-real}
In case \xref{cases-II-q5-B2}, $X\simeq \PP(1,1,1,2)$.
\end{tparag}
\begin{proof}
Let $S\in |2L|$ be a general member.
Then $S$ is Cartier and by Proposition \ref{prop-base-point-H}
$X$ is has at worst Du Val singularities.
By the adjunction formula $S$ is a del Pezzo surface of 
degree $9$. It follows that
$S$ is smooth and $S\simeq \PP^2$ (see Remark \ref{rem-class-delPezzo}).
The restriction map $H^0(X,\OOO_X(S))\to H^0(S,\OOO_S(S))$
is surjective. Hence the linear system $|S|$ is base point free
and determines a morphism $\varphi\colon X\to \PP^6$.
We have $(\deg \varphi)(\deg \varphi (X))=S^3=4$.
So $\varphi$ is birational and $\varphi(X)\subset \PP^6$
is a variety of degree $4$.
A general hyperplane section $\varphi(S)\subset \varphi(X)$
is a Veronese surface. It is well-known that in this situation
$\varphi(X)$ is a cone over $\varphi(S)$, i.e.,
$X\simeq \varphi(X)\simeq \PP(1,1,1,2)$.
\end{proof}
\end{mparag}

\begin{mparag}{Cases \ref{cases-I-3-2},
\ref{cases-II-q4-B5},
\ref{cases-II-q4-B55},
\ref{cases-II-q5-B26},
\ref{cases-II-q5-B7},
\ref{cases-II-q5-B2236},
\ref{cases-II-q7-B39}, 
\ref{cases-II-q7-B210}.}
\label{case-last-mnogo}
We apply construction \ref{construction-1a}.
Let $r$ be the Gorenstein index of $X$.
First we construct a birational
extremal extraction $g\colon \tilde X\to X$
such that $\tilde X$ has only terminal singularities and
the exceptional divisor $E$ of $g$ has discrepancy $1/r$. 

\begin{tparag}{Claim.}
\label{claim-blowup}
Either 
\begin{enumerate}
\item 
There is a cyclic quotient singularity $P\in X$ 
of type $\frac 1r (b,-b,1)$, where $\gcd (r,b)=1$, or
\item
we are in case \xref{cases-II-q4-B55}
and there is a point $P\in X$ 
of type $cA/5$ of axial weight $2$.
\end{enumerate}
\end{tparag}
\begin{proof}
Note that in all cases there is a basket 
point $P\in \B$
of index $r$.
If this point is unique, it corresponds to a cyclic quotient singularity
of $X$.
The point $P\in \B$ of index $r$ is not unique only in case
\ref{cases-II-q4-B55}. Then $r=5$ and there are two 
points $P_1,\, P_2\in \B$ of index $5$.
They correspond either two cyclic quotient singularities of $X$ 
or a point $P\in X$ 
of type $cA/5$. 
\end{proof}
In case (i) the weighted blowup 
of $P\in X$ with weights $\frac 1r (b,r-b,1)$
gives us a desired contraction $g$. Similarly, in case (ii)
a suitable weighted blowup 
gives us a desired contraction $g$ (see \cite{Kawamata-1992-e-app}).

Further, $r\HHH$ is the linear system of Cartier divisors.
Hence we can write $g^*\HHH=\tilde \HHH+\delta E$, where $\delta\ge 1/r$.
Thus,
\begin{equation}
\label{eq-beliberda-crepant}
-K_{\tilde X}\qq g^*(-K_X)-{\textstyle{\frac1r}}E\qq \tilde \HHH+(\delta-{\textstyle{\frac1r}}) E.
\end{equation}
By Corollary \ref{cor--K-isol-base-points}
the linear system $\tilde \HHH$ has only isolated base points outside of $E$.
Therefore, $-K_{\tilde X}$ is nef.

If $g(E)$ is a cyclic quotient singularity,
then $E\simeq \PP(b,r-b,1)$, $E|_E\sim \OOO_{\PP(b,r-b,1)}(-r)$, 
and $E^3=r^2/b(r-b)$. Therefore,
\[
-K_{\tilde X}^3=-K_X^3-{\textstyle{\frac1{r^3}}}E^3\ge 
\frac{125}2-\frac {r^2}{b(r-b)}>0.
\]
This shows that $-K_{\tilde X}$ is big.
Similar computations shows that this fact also holds in case
\ref{claim-blowup}, (ii).

Let $C$ be a curve such that $-K_{\tilde X}\cdot C=0$.
By \eqref{eq-beliberda-3} we have 
$q\tilde L\cdot C+\beta E\cdot C=0$.
By \eqref{eq-beliberda-crepant} $E\cdot C> 0$.
Hence $\tilde L\cdot C<0$. Since $\dim |L|>0$,
there is at most a finite number of such curves.
Thus the linear system $|-nK_{\tilde X}|$ does not contract any divisors.

\begin{parag}{}
\label{eq-beliberda-tabl}
Consider diagram \eqref{eq-main-diag-1}.
Since $K_X+qL\sim 0$, the constant $\beta$ in \ref{eq-beliberda-beta}
is a non-negative integer.
We can write
\[
K_{\tilde X}=g^*K_X+\frac 1r E, \qquad
\tilde L=g^*L-\delta E,
\]
where $\delta\in \QQ$, $\delta>0$.
Since $rL$ is Cartier (see Lemma \ref{lemma-Kawamata-index-divisors}),
$\delta=k/r$ for some $k\in \ZZ$, $k>0$. Therefore,
\[
\beta =-\frac 1r+q\delta=\frac{qk-1}r
\]
and the value of $\beta$ 
is bounded from below as follows:
\begin{equation*}
\begin{array}{|c||c|c|c|c|}
\hline
&&&&
\\
\text{case}
&\text{\ref{cases-I-3-2}} 
&\text{\ref{cases-II-q4-B5}}
\text{\ref{cases-II-q4-B55}}
\text{\ref{cases-II-q7-B39}}
&\text{\ref{cases-II-q5-B26}}
\text{\ref{cases-II-q5-B2236}}
&\text{\ref{cases-II-q5-B7}}
\text{\ref{cases-II-q7-B210}}
\\[9pt]
\hline
&&&&
\\
\beta &\ge 1&\ge 3&\ge 4&\ge 2
\\[9pt]
\hline
\end{array}
\end{equation*}
\end{parag}

\begin{parag}{}
\label{par-birat-case-first}
First we assume that $\dim Z=\dim X$.
Then $f$ is a divisorial contraction and $Z$ 
is a $\QQ$-Fano threefold. 
By \eqref{eq-beliberda-3}
we have $K_Z+qL_Z+\beta E_Z\qq 0$, where $E_Z$ and $L_Z$ 
are effective non-zero divisors. Hence, 
$q\QQ(Z)\ge q+\beta >4$. 
In particular, $Z$ is not Gorenstein (see 
Corollary \ref{cor-first-prop-index}).

Assume that $-K_Z^3< 125/2$.
By Proposition \ref{prop-diag-0}
$\dim |-K_X|=\dim |-K_Z|=32$.
Hence $X$ is of type \ref{cases-II-q5-B2236}.
By \eqref{eq-KZ-K3} $\dim |-K_Z|\ge 60$
and by \eqref{eq-beliberda-3} $q\QQ(Z)\ge 9$.
On the other hand,
$\discrep(Z)\ge \discrep(\tilde X)\ge 1/5$.
Therefore the Gorenstein index of $Z$
is at most $5$ (see \cite{Kawamata-1992-e-app}).
By Proposition \ref{proposition-Cl-2} $\Cl Z\simeq \ZZ$.
Let $L'$ be the ample generator of $\Cl Z\simeq \ZZ$,
let $r'\le 5$ be the Gorenstein index of $Z$,
and let $S\in |r' L'|$ a general member.
Then $S$ be the ample generator 
of $\Pic Z$. By Proposition \ref{prop-base-point-H}
$S$ has at worst Du Val singularities.
By the adjunction formula
$K_S=(r'-q\QQ(Z)) L'|_S$.
Since $L'|_S$ is a Cartier divisor, $S$ is 
a del Pezzo surface with $qF(S)\ge q\QQ(Z)-r'\ge 4$.
This is impossible (see \ref{rem-class-delPezzo}).
Thus $-K_Z^3\ge 125/2$ and $Z$ is such as in \ref{cases-II}.

Now we consider possibilities for $X$ case by case.
In cases
\ref{cases-II-q5-B26},
\ref{cases-II-q5-B2236},
\ref{cases-II-q7-B39}, and
\ref{cases-II-q7-B210} we have $q\QQ(Z)\ge 9$, a contradiction.
In cases \ref{cases-II-q4-B5},
\ref{cases-II-q4-B55}, and
\ref{cases-II-q5-B7} we have $q\QQ(Z)= 7$. 
Hence $Z$ is such as in \ref{cases-II-q7-B39} or
\ref{cases-II-q7-B210}.
Then $q+\beta=7$. By 
\eqref{eq-beliberda-3} $L_Z$ and $E_Z$ are linear equivalent and 
they are generators of $\Cl Z$. 
On the other hand, $\dim |L| \ge 2>\dim |L_Z|=1$, a contradiction.

In case \ref{cases-I-3-2} 
$\tilde X$ is of Gorenstein index $2$. Hence,
$\discrep(\tilde X)=1/2$.
On the other hand, 
$f{\comp}\chi$
is a composition of a flop and steps of the $K$-MMP.
Therefore, $\discrep(Z)\ge 1/2$. 
This is possible only if $Z$ of type \ref{cases-II-q5-B2}.
But then $35=\dim |-K_X|>\dim |-K_Z|=33$,
a contradiction.
\end{parag}

\begin{parag}{}
Thus we may assume that $\dim Z<\dim X$.
Let $M\in |2L|$ be a general member.
Note that by \ref{eq-beliberda-tabl} $q+\beta \ge 3$
and $q+\beta =3$ only in case \ref{cases-I-3-2}.
By \eqref{eq-beliberda-4}
$L_{Y}$ can be $f$-horizontal only in case \ref{cases-I-3-2}
and if $Z$ is a curve.
By Lemma \ref{cor-dim1-index3} we have a contradiction.
Hence $L_{Y}$ is $f$-vertical.
As in Lemma \ref{lemma-constr-diag-descr} we have
$L_{Y}=f^*\Xi$ for 
some integral Weil divisor $\Xi$ on $Z$,
$\dim |\Xi|=\dim |L|\ge 1$, and $\Xi$ is a generator of $\Cl Z/\mathrm{Tors}$.
\end{parag}

\begin{parag}{}
Assume that $Z$ is a surface. 
From \eqref{eq-beliberda-4} we get $\beta \le 2$.
By \ref{eq-beliberda-tabl} this is possible only in cases 
\ref{cases-I-3-2}, \ref{cases-II-q5-B7} or \ref{cases-II-q7-B210}.
If $K_Z^2<8$, we have 
$\dim |\Xi|=0$, a contradiction.
Hence $Z$ is either $\PP^2$ or $\PP(1,1,2)$.
Consider the case $Z\simeq \PP(1,1,2)$.
Then $\dim |\Xi|=1$ and we are in case
\ref{cases-II-q7-B210}. 
Let $M\in |3L|$ be a general member.
We can write 
$K_{Y}+2M_{Y}+L_{Y}+\gamma E_{Y}\sim 0$, where $\gamma>0$.
This shows that $M_{Y}$ is $f$-vertical.
Thus $M_{Y}\sim 3L_{Y}=3f^*\Xi$
and $\dim |M_{Y}|=\dim |3\Xi|=4$, a contradiction.

Consider the case $Z\simeq \PP^2$.
Then $\dim |\Xi|=2$ and we are in case \ref{cases-II-q5-B7}. 
Let $M\in |2L|$ be a general member.
We can write 
$K_{Y}+2M_{Y}+L_{Y}+\gamma E_{Y}\sim 0$, where $\gamma>0$.
This shows that $\gamma=\beta=2$ and $M_{Y}$ is $f$-vertical.
Thus $M_{Y}\sim 2L_{Y}=f^*\Xi$
and $\dim |M_{Y}|=\dim |2\Xi|=5$, a contradiction.
\end{parag}

\begin{parag}{}
Assume that $Z$ is a curve.
Then $Z\simeq \PP^1$. Since $L_Y=f^*\Xi$ is not divisible 
in $\Cl Y$, $\dim |\Xi|\le 1$. So we are in cases 
\ref{cases-II-q5-B2236},
\ref{cases-II-q7-B39}, or
\ref{cases-II-q7-B210}.
Moreover, since $\dim |L|>0$, $\dim |\Xi|= 1$.
Case \ref{cases-II-q5-B2236} is impossible because then $\beta \ge 4$.
Let $M\in |2L|$ be a general member.
We can write 
$K_{Y}+3M_{Y}+L_{Y}+\gamma E_{Y}\sim 0$, where $\gamma>0$.
This shows that $M_{Y}$ is $f$-vertical.
Thus $M_{Y}\sim 2L_{Y}=2f^*\Xi$
and $\dim |M_{Y}|=\dim |2\Xi|=2$, a contradiction.
\end{parag}
\end{mparag}

Now we consider case \ref{cases-II-q6-B57}.
\begin{mparag}{Case \ref{cases-II-q6-B57}.}
\label{case-last-5.1.7}
By Lemmas \ref{lemma-canonical-2} and \ref{lemma-canonical-3}
the pair $(X,|-K_X|)$ is not canonical. Thus we apply the construction 
\eqref{eq-main-diag-1} in case \ref{construction-1b}.
Then in \eqref{eq-crepant-formula} we have $\alpha>0$.
Assume that $\dim Z=3$. Since $\alpha>0$, and by 
Proposition \ref{prop-A-lin-syst} we have
$\dim |-K_Z|>\dim |-K_X|=38$. Then 
by Proposition \ref{prop-diag-0} 
$-K_Z^3\ge 125/2$. 
Hence $Z$ is $\QQ$-Fano 
such as in Proposition \ref{cases-II}.
Moreover, by \eqref{eq-beliberda-3} we have
$q\QQ(Z)\ge q\QQ(X)+\beta=6+\beta$.
This implies that 
$E_Z\sim L_Z$ is a generator of $\Cl Z$, $q\QQ(Z)=7$, and $\beta=1$.
So,
the variety $Z$ is of type \ref{cases-II-q7-B39}.
Obviously, $\dim |2L_Z|\ge \dim |2L|$.
This contradicts Proposition \ref{cases-II}.

Thus $\dim Z=1$ or $2$.
If $Z$ is a surface, then by
Lemma \ref{lemma-constr-diag-descr} 
$Z\simeq\PP(1,1,2)$.
Let $M\in |2L|$ be a general member.
We can write 
$K_{Y}+3M_{Y}+\gamma E_{Y}\sim 0$, where $\gamma>0$.
Restricting to a general fibre
we obtain that $M_{Y}$ is $f$-vertical.
Thus, $M_{Y}\sim 2L_{Y}=2f^*\Xi$
and $\dim |M_{Y}|=\dim |2\Xi|\le 3$, a contradiction.
\end{mparag}

Finally we consider cases when $q\QQ(X)=1$.
\begin{mparag}{Cases \ref{cases-I-2}, \ref{cases-I-2-2}, \ref{cases-I-3-1}.}
\label{case-last-index1}
By Lemmas \ref{lemma-canonical-2} and \ref{lemma-canonical-3} the pair
$(X,|-K_X|)$ is not canonical. Thus we may apply 
construction \ref{assumption-Sarkisov-link} under assumptions 
\ref{construction-1b}.

Then in \eqref{eq-crepant-formula} we have $\alpha>0$.
Assume that $\dim Z=3$. 
Similar to \ref{case-last-5.1.7}
$\dim |-K_Z|>\dim |-K_X|$
and
$-K_Z^3\ge 125/2$. 
Hence $Z$ is $\QQ$-Fano 
such as in Proposition \ref{cases-II} or \ref{cases-I} with
$q\QQ (Z)>1$.
By \ref{subsect-case-last-real},
\ref{case-last-mnogo}, and \ref{case-last-5.1.7} 
$Z$ is of type \ref{cases-II-q5-B2} and $Z\simeq \PP(1,1,1,2)$.
Then $\dim |-K_X|<\dim |-K_Z|=33$,
so $X$ is of type \ref{cases-I-2} and $\dim \HHH_Z\ge 32$.
Easy computations show that $\HHH_Z\sim \OOO_{\PP(1,1,1,2)}(n)$,
with $n\ge 5$. 
On the other hand, $-K_Z\sim \HHH_Z +\alpha E_Z$,
where $\alpha>0$, a contradiction.

Therefore,
$1\le \dim Z\le 2$. If $Z$ is a curve,
we have a contradiction by 
Lemma \ref{cor-dim1-index3}. 
Thus $Z$ is a surface. 
Then by Lemma \ref{lemma-constr-diag-dimZ2-not-canon}
the fibration $f$ has no discriminant curve.
Hence by Lemma \xref{lemma-constr-diag-dimZ2-q1}
we have $\dim |-K_X|<30$, a contradiction.
\end{mparag}


\begin{thebibliography}{10}

\bibitem{Kawamata-1992bF}
Y. Kawamata.
\newblock Boundedness of {$\bold Q$}-{F}ano threefolds.
\newblock In {\em Proceedings of the International Conference on Algebra, Part
  3 (Novosibirsk, 1989)}, volume 131 of {\em Contemp. Math.}, pages 439--445,
  Providence, RI, 1992. Amer. Math. Soc.

\bibitem{Namikawa-1997}
Y. Namikawa.
\newblock Smoothing {F}ano {$3$}-folds.
\newblock {\em J. Algebraic Geom.}, 6(2):307--324, 1997.

\bibitem{Reid-YPG1987}
M. Reid.
\newblock Young person's guide to canonical singularities.
\newblock In {\em Algebraic geometry, Bowdoin, 1985 (Brunswick, Maine, 1985)},
  volume~46 of {\em Proc. Sympos. Pure Math.}, pages 345--414. Amer. Math.
  Soc., Providence, RI, 1987.

\bibitem{Suzuki-2004}
K. Suzuki.
\newblock On {F}ano indices of {$\Bbb Q$}-{F}ano 3-folds.
\newblock {\em Manuscripta Math.}, 114(2):229--246, 2004.

\bibitem{Alexeev-1994ge}
V. Alexeev.
\newblock General elephants of {${\bf Q}$}-{F}ano 3-folds.
\newblock {\em Compositio Math.}, 91(1):91--116, 1994.

\bibitem{Takagi-2002-I-II}
H. Takagi.
\newblock On classification of {$\Bbb Q$}-{F}ano 3-folds of {G}orenstein index
  2. {I}, {II}.
\newblock {\em Nagoya Math. J.}, 167:117--155, 157--216, 2002.

\bibitem{Takagi-2006}
H. Takagi.
\newblock Classification of primary {$\Bbb Q$}-{F}ano threefolds with
  anti-canonical {D}u {V}al {$K3$} surfaces. {I}.
\newblock {\em J. Algebraic Geom.}, 15(1):31--85, 2006.

\bibitem{Prokhorov-2005a}
Yu. Prokhorov.
\newblock On the degree of {F}ano threefolds with canonical {G}orenstein
  singularities.
\newblock {\em Russian Acad. Sci. Sb. Math.}, 196(1):81--122, 2005.

\bibitem{Prokhorov-2006-Enr}
Yu. Prokhorov.
\newblock On {F}ano-{E}nriques threefolds, 2006.

\bibitem{Iskovskikh-Prokhorov-1999}
V.~A. Iskovskikh and Yu.~G. Prokhorov.
\newblock {\em Fano varieties. {A}lgebraic geometry. {V}.}, volume~47 of {\em
  Encyclopaedia Math. Sci.}
\newblock Springer, Berlin, 1999.

\bibitem{Kawamata-1988-crep}
Y. Kawamata.
\newblock Crepant blowing-up of {$3$}-dimensional canonical singularities and
  its application to degenerations of surfaces.
\newblock {\em Ann. of Math. (2)}, 127(1):93--163, 1988.

\bibitem{KMMT-2000}
J. Koll{\'a}r, Yoichi Miyaoka, Shigefumi Mori, and Hiromichi Takagi.
\newblock Boundedness of canonical {$\bold Q$}-{F}ano 3-folds.
\newblock {\em Proc. Japan Acad. Ser. A Math. Sci.}, 76(5):73--77, 2000.

\bibitem{Mori-Prokhorov-2006}
S.~Mori and Yu. Prokhorov.
\newblock On $\mathbf {Q}$-conic bundles, 2006.

\bibitem{Fujita-1975}
T. Fujita.
\newblock On the structure of polarized varieties with {$\Delta $}-genera zero.
\newblock {\em J. Fac. Sci. Univ. Tokyo Sect. IA Math.}, 22:103--115, 1975.

\bibitem{Mukai-2002}
Sh. Mukai.
\newblock New developments in the theory of {F}ano threefolds: vector bundle
  method and moduli problems [translation of {S}\=ugaku {\bf 47} (1995), no.\
  2, 125--144].
\newblock {\em Sugaku Expositions}, 15(2):125--150, 2002.

\bibitem{Ambro-1999}
F.~Ambro.
\newblock Ladders on {F}ano varieties.
\newblock {\em J. Math. Sci. (New York)}, 94(1):1126--1135, 1999.
\newblock Algebraic geometry, 9.

\bibitem{Kollar-ShB-1988}
J.~Koll{\'a}r and N.~I. Shepherd-Barron.
\newblock Threefolds and deformations of surface singularities.
\newblock {\em Invent. Math.}, 91(2):299--338, 1988.

\bibitem{Kawachi-Masek-1998}
T. Kawachi and V. Ma{\c{s}}ek.
\newblock Reider-type theorems on normal surfaces.
\newblock {\em J. Algebraic Geom.}, 7(2):239--249, 1998.

\bibitem{Reider-1988}
I. Reider.
\newblock Vector bundles of rank {$2$} and linear systems on algebraic
  surfaces.
\newblock {\em Ann. of Math. (2)}, 127(2):309--316, 1988.

\bibitem{Sano-1996}
T. Sano.
\newblock Classification of non-{G}orenstein {${\bf Q}$}-{F}ano {$d$}-folds of
  {F}ano index greater than {$d-2$}.
\newblock {\em Nagoya Math. J.}, 142:133--143, 1996.

\bibitem{Kawamata-1992-e-app}
Y. Kawamata.
\newblock The minimal discrepancy coefficients of terminal singularities in
  dimension three ({A}ppendix to {V}{.}{V}{.} {S}hokurov{'}s paper {"}3-fold
  log flips{"}).
\newblock {\em Russ. Acad. Sci., Izv., Math.}, 40(1):193--195, 1993.

\end{thebibliography}

\def\cprime{$'$}

\end{document}